\input amstex  
\documentstyle{amsppt}
\magnification=\magstephalf  
 \addto\tenpoint{\baselineskip 15pt
  \abovedisplayskip18pt plus4.5pt minus9pt
  \belowdisplayskip\abovedisplayskip
  \abovedisplayshortskip0pt plus4.5pt
  \belowdisplayshortskip10.5pt plus4.5pt minus6pt}\tenpoint
\pagewidth{6.5truein} \pageheight{8.9truein}
\subheadskip\bigskipamount
\belowheadskip\bigskipamount
\aboveheadskip=3\bigskipamount
\catcode`\@=11
\def\output@{\shipout\vbox{%
 \ifrunheads@ \makeheadline \pagebody
       \else \pagebody \fi \makefootline 
 }%
 \advancepageno \ifnum\outputpenalty>-\@MM\else\dosupereject\fi}
\outer\def\subhead#1\endsubhead{\par\penaltyandskip@{-100}\subheadskip
  \noindent{\subheadfont@\ignorespaces#1\unskip\endgraf}\removelastskip
  \nobreak\medskip\noindent}
\outer\def\enddocument{\par
  \add@missing\endRefs
  \add@missing\endroster \add@missing\endproclaim
  \add@missing\enddefinition
  \add@missing\enddemo \add@missing\endremark \add@missing\endexample
 \ifmonograph@ 
 \else
 \vfill
 \nobreak
 \thetranslator@
 \count@\z@ \loop\ifnum\count@<\addresscount@\advance\count@\@ne
 \csname address\number\count@\endcsname
 \csname email\number\count@\endcsname
 \repeat
\fi
 \supereject\end}
\catcode`\@=\active
\CenteredTagsOnSplits
\NoBlackBoxes
\nologo
\def\today{\ifcase\month\or
 January\or February\or March\or April\or May\or June\or
 July\or August\or September\or October\or November\or December\fi
 \space\number\day, \number\year}
\define\({\left(}
\define\){\right)}

\define\Aut{\operatorname{Aut}}
\define\CC{{\Bbb C}}
\define\CP{{\Bbb C\Bbb P}}

\define\ZZ{{\Bbb Z}}
\define\[{\left[}
\define\]{\right]}

\define\chiup{\raise.5ex\hbox{$\chi$}}

\define\exertag #1#2{#2\ #1}

\define\ind{\operatorname{ind}}
\define\inv{^{-1}}
\define\mstrut{^{\vphantom{1*\prime y}}}
\define\protag#1 #2{#2\ #1}

\define\res#1{\negmedspace\bigm|_{#1}}
\define\temsquare{\raise3.5pt\hbox{\boxed{ }}}

\define\theprotag#1 #2{#2~#1}

\define\xca#1{\removelastskip\medskip\noindent{\smc%
#1\unskip.}\enspace\ignorespaces }

\define\zmod#1{\ZZ/#1\ZZ}

\define\zt{\zmod2}

\NoRunningHeads 

\define\Fred{\operatorname{Fred}}
\define\Hh{\check H}
\define\SAFred{\operatorname{Fred}^1}
\define\Sym{\operatorname{Sym}}
\define\TT{\Bbb{T}}
\define\Vect{\operatorname{Vect}}
\define\bul{\bullet}
\define\lh{\check \lambda }
\define\sE{\Cal{E}}
\define\sT{\Cal{T}}

\refstyle{A}
\widestnumber\key{SSSSS}   
\document

	\topmatter
 \title\nofrills The Verlinde Algebra is Twisted Equivariant $K$-Theory
\endtitle 
 \author Daniel S. Freed  \endauthor
 \thanks The author is supported by NSF grant DMS-0072675.\endthanks
 \affil Department of Mathematics \\ University of Texas at Austin\endaffil 
 \address Department of Mathematics, University of Texas, Austin, TX
78712\endaddress 
 \email dafr\@math.utexas.edu \endemail
 \date January 4, 2001\enddate
	\endtopmatter

\document

 \comment
 lasteqno 0@  1
 \endcomment

$K$-theory in various forms has recently received much attention in
10-dimensional superstring theory.  Our raised consciousness about {\it
twisted\/} $K$-theory led to the serendipitous discovery that it enters in a
different way into 3-dimensional topological field theories, in particular
Chern-Simons theory.  Namely, as the title of the paper reports, the Verlinde
algebra is a certain twisted $K$-theory group.  This assertion, and its
proof, is joint work with Michael Hopkins and Constantin Teleman.  The
general theorem and proof will be presented elsewhere~\cite{FHT}; our goal
here is to explain some background, demonstrate the theorem in a simple
nontrivial case, and motivate it through the connection with topological
field theory.
 
From a mathematical point of view the Verlinde algebra is defined in the
theory of loop groups.  Let $G$~be a compact Lie group.  There is a version
of the theorem for {\it any\/} compact group~$G$, but here for the most part
we focus on connected, simply connected, and simple groups---$G=SU_2$ is the
simplest example.  In this case a central extension of the free loop
group~$LG$ is determined by the {\it level\/}, which is a positive
integer~$k$.  There is a finite set of equivalence classes of positive energy
representations of this central extension; let $V_k(G)$~denote the free
abelian group they generate.  One of the influences of 2-dimensional
conformal field theory on the theory of loop groups is the construction of an
algebra structure on~$V_k(G)$, the {\it fusion product\/}.  This is the {\it
Verlinde algebra\/}~\cite{V}.
 
Let $G$~act on itself by conjugation.  Then with our assumptions the
equivariant cohomology group~$H^3_G(G)$ is free of rank one.  Let $h(G)$~be
the dual Coxeter number of~$G$, and define~ $\zeta (k)\in H^3_G(G)$ to be
$k+h(G)$~times a generator.  We will see that elements of~$H^3$ may be used
to {\it twist\/} $K$-theory, and so elements of equivariant~$H^3$ twist
equivariant $K$-theory.

        \proclaim{\protag{(Freed-Hopkins-Teleman)} {Theorem}}
 There is an isomorphism of algebras 
  $$ V_k(G)\cong K^{\dim G + \zeta (k)}_{G}(G),  $$
where the right hand side is the $\zeta (k)$-twisted equivariant $K$-theory
in degree~$\dim G$.
        \endproclaim
 
\flushpar
 The group structure on the right-hand side is induced from the
multiplication map $G\times G\to G$.  

For an arbitrary compact Lie group~$G$ the level~$k$ is replaced by a class
in~$H^4(BG;\ZZ)$ and the dual Coxeter number~$h(G)$ is pulled back from a
universal class in~$H^4(BSO;\ZZ)$ via the adjoint
representation.\footnote{The story is a bit more subtle and
involves~$H^2(BG;\zt)$---so a twisting using~$H^1_G(G;\zt)$---as well,
though we do not discuss that here.}  The twisting class is obtained from
their sum by transgression.

I warmly thank Mike Hopkins and Constantin Teleman for their continued
collaboration and for comments on this manuscript.

 \newpage
 \head
 \S{1} Twisted $K$-theory
 \endhead
 \comment
 lasteqno 1@ 11
 \endcomment

Twistings of cohomology theories are most familiar for ordinary cohomology.
Let $M$~be a manifold (or suitably nice space).  Then a flat real vector
bundle $E\to M$ determines twisted real cohomology groups~$H^{\bul}(M;E)$.
In differential geometry these cohomology groups are defined by extending the
de Rham complex to forms with coefficients in~$E$ using the flat connection.
The sorts of twistings of $K$-theory we consider are one-dimensional, so
analogous to the case when $E$~is a line bundle.  There are also
one-dimensional twistings of {\it integral\/} cohomology, determined by a
{\it local system\/} $Z\to M$.  This is a bundle of groups isomorphic
to~$\ZZ$, so is determined up to isomorphism by an element
of~$H^1\bigl(M;\Aut(\ZZ)\bigr)\cong H^1(M;\zt)$, since the only nontrivial
automorphism of~$\ZZ$ is multiplication by~$-1$.  The twisted integral
cohomology $H^\bul(M;Z)$ may be thought of as sheaf cohomology, or defined
using a cochain complex.  We give a \v Cech description as follows.  Let
$\{U_i\}$ be an open covering of~$M$ and
  $$ g_{ij}\: U_i\cap U_j\longrightarrow  \{\pm1\} \tag{1.1} $$
a cocycle defining the local system~$Z$.  Then an element of~$H^q(M;Z)$ is
represented by a collection of $q$-cochains $a_i\in Z^q(U_i)$ which satisfy 
  $$ a_j = g_{ij}\,a_i\qquad \text{on $U_{ij}=U_i\cap U_j$}. \tag{1.2} $$
We can use any model of cochains, since the group ~$\Aut(\ZZ)\cong \{\pm1\}$
always acts.  In place of cochains we represent integral cohomology classes
by maps to an Eilenberg-MacLane space~$K(\ZZ,q)$.  The cohomology group is
the set of homotopy classes of maps, but here we use honest maps as
representatives.  The group~$\Aut(\ZZ)$ acts on~$K(\ZZ,q)$.  One model
of~$K(\ZZ,0)$ is the integers, with $-1$~acting by multiplication.  The
circle is a model for~$K(\ZZ,1)$, and $-1$~acts by reflection.  Using the
action of~$\Aut(\ZZ)$ on~$K(\ZZ,q)$ and the cocycle~\thetag{1.1} we build an
associated bundle $\Cal{H}^q \to M$ with fiber~$K(\ZZ,q)$.
Equation~\thetag{1.2} says that twisted cohomology classes are represented by
sections of~$\Cal{H}^q\to M$; the twisted cohomology group~$H^q(M;Z)$ is the
set of homotopy classes of sections of~$\Cal{H}^q\to M$.
 
Twistings may be defined for any generalized cohomology theory; our interest
is in complex $K$-theory.  In homotopy theory one regards $K$ as a marriage
of a ring and a space (more precisely, spectrum), and it makes sense to ask
for the units in~$K$, denoted~$GL_1(K)$.  In the previous paragraph we used
the units in integral cohomology, the group~$\zt$.  For complex $K$-theory
there is a richer group of units\footnote{This is a consequence of results
in~\cite{DK}, \cite{S}, \cite{ASe}, \cite{AP}.}
  $$ GL_1(K) \sim \zt \times \CP^{\infty} \times BSU. \tag{1.3} $$
In our problem the last factor doesn't enter and all the interest is in the
first two, which we denote~$GL_1(K)'$.  As a first approximation, view~$K$ as
the category of all finite dimensional $\zt$-graded complex vector spaces.
Then $\CP^{\infty}$~is the subcategory of even complex lines, and it is a
group under tensor product.  It acts on~$K$ by tensor product as well.  The
nontrivial element of~$\zt$ in~\thetag{1.3} acts on~$K$ by reversing the
parity of the grading.  This model is deficient since there is not an
appropriate topology.  One may consider instead complexes of complex vector
spaces, or spaces of operators as we do below.  Of course, there are good
topological models of~$\CP^{\infty}$, for example the space of all
one-dimensional subspaces of a fixed complex Hilbert space.  For a
manifold~$M$ the twistings of $K$-theory of interest are classified up to
isomorphism by
  $$ H^1\bigl(M;GL_1(K)' \bigr) \cong H^1(M;\zt)\times H^3(M;\ZZ). \tag{1.4}
     $$
In this paper we will not encounter twistings from the first factor and will
focus exclusively on the second.\footnote{In passing we remark that twistings
form an abelian group, and the isomorphism~\thetag{1.4} is {\it not\/} an
isomorphism of groups; rather the sum of twisting classes in~$H^1(M;\zt)$ has
a component in~$H^3(M;\ZZ)$.}  These twistings are represented by
cocycles~$g_{ij}$ with values in the space of lines, in other words by
complex line bundles $L_{ij}\to U_{ij}$ which satisfy a cocycle condition.
This is the data often given to define a {\it gerbe\/}.

Following~\cite{A}, we present a model of twisted $K$-theory in terms of
operators on an infinite dimensional separable complex Hilbert space~$H$.
Let $PGL(H)$ be the projective general linear group of Hilbert space.
Kuiper's theorem asserts that the set of invertible transformations~$GL(H)$
is contractible, whence $PGL(H)= GL(H)/\CC^{\times }$ has the homotopy type
of~$\CP^{\infty}$.  The space of Fredholm operators~$\Fred^0(H)$ has the
homotopy type of $\ZZ\times BU$, so is a classifying space for~$K^0$, and a
suitable space of self-adjoint Fredholm operators~$\SAFred(H)$ is a
classifying space for~$K^1$ (see~\cite{AS}).  (Alternatively, we could take
invertible operators of the form $1 + \text{compact}$ as a classifying space
for~$K^1$.)  Recall Bott periodicity which implies that $K^q$~depends only on
the parity of~$q$, so we need only consider~$q=0$ and~$q=1$.  Then a twisting
class in~$H^3(M;\ZZ)$ is represented by a cocycle
  $$ g_{ij}\: U_{ij}\longrightarrow PGL(H),  $$
and a twisted $K$-theory class by maps 
  $$ a_i \: U_i\longrightarrow \Fred^q(H)  $$
which satisfy 
  $$ a_j = g_{ij}a_ig_{ij}\inv \qquad \text{on $U_{ij}$}.  $$
Associated to~$\{g_{ij}\}$ is a bundle over~$M$ with fiber~$\Fred^q(H)$, and
the twisted $K$-groups are the sets of homotopy classes of sections.
 
It is important that twisted cohomology theories (such as twisted $K$-theory)
are cohomology theories.  In particular, they satisfy the Mayer-Vietoris
property.  Note that twisted cohomology is not a ring; rather, it is a module
over untwisted cohomology.  When twisted cohomology classes are multiplied
the twisting class adds.
 
There are also twistings of equivariant cohomology theories.  Suppose a
compact Lie group~$G$ acts on~$M$.  Then the twistings of equivariant
$K$-theory we consider here are classified by the equivariant cohomology
group~$H^3_G(M;\ZZ)$.  

We conclude this section by computing the twisted equivariant $K$-theory
of~$SU_2$, where $SU_2$~acts on itself by conjugation.  As a preliminary, we
first recall some facts about equivariant $K$-theory.  Let $R(G)$~denote the
representation ring of~$G$.  Then
  $$ K^q_G(pt) \cong \cases R(G) ,&q=0;\\0,&q=1.\endcases \tag{1.5} $$
Now suppose $H\subset G$ is a closed subgroup.  Then 
  $$ K_G(G/H) \cong K_H(pt). \tag{1.6} $$
The isomorphism is restriction to the basepoint; its inverse is the
associated bundle construction.  The computation below is somewhat easier in
$K$-homology, rather than $K$-cohomology.  Poincar\'e duality gives
properties of equivariant $K$-homology analogous to~\thetag{1.5}
and~\thetag{1.6}.
 
The representation ring of~$SU_2$ is a polynomial ring generated by the
two-dimensional representation~$\sigma $.  Let $\TT\subset SU_2$ denote the
maximal torus of diagonal matrices.  Then $R(\TT)\cong \ZZ[\alpha ,\alpha
\inv ]$, where $\alpha $~is the one-dimensional standard representation.
There is a holomorphic induction map 
  $$ \ind\:R(\TT)\longrightarrow R(SU_2), \tag{1.7} $$
also known as the Borel-Weil construction.  Namely, a representation of~$\TT$
determines a holomorphic line bundle over~$SU_2/\TT\cong \CP^1$, and the
induced representation of~$SU_2$ is the space of its holomorphic sections.
For example, $\ind(\alpha )=\sigma $, $\ind(1) = 1$, and $\ind(\alpha \inv
)=0$.  Now $R(\TT)$ is an $R(SU_2)$-module---a representation of~$SU_2$
multiplies a representation of~$\TT$ by restriction---which is free of
rank~2.  We take $1,\alpha \inv $ as generators.  Note that
\thetag{1.7}~preserves the $R(SU_2)$-module structure.
 
We are ready to compute.  As stated earlier we let $G=SU_2$~act on itself by
conjugation.  Let $\zeta \in H^3_{G}(G)\cong \ZZ$ be $m$~times the
generator.  Cover~$SU_2$ by the invariant sets 
  $$ \aligned
      U &= SU_2\setminus \{-1\} \\ 
      V &= SU_2\setminus \{+1\}.\endaligned  $$
Then  
  $$ \aligned
      U, V &\sim pt \\ 
      U\cap V &\sim SU_2/\TT.\endaligned \tag{1.8} $$
The twisting class is represented by the $m^{\text{th}}$~power of the
hyperplane bundle on~$U\cap V\sim\CP^1$, which is induced from the
representation~$\alpha ^m$ of~$\TT$.  Let $K^G_{q+\zeta }(G)$ denote the
$\zeta $-twisted equivariant $K$-homology in degree~$q$.  Using~\thetag{1.5}
and~\thetag{1.8} we write the Mayer Vietoris sequence for homology as
  $$ 0 \longrightarrow K^G_{1+\zeta }(G) \longrightarrow K^G_0(U\cap V)
     \longrightarrow K^G_0(U) \oplus K^G_0(V) \longrightarrow K^G_{0+\zeta
     }(G) \longrightarrow 0. \tag{1.9} $$
The middle terms are untwisted since the twisting class~$\zeta $ restricts
trivially on~$U,V$.  From~\thetag{1.8}, \thetag{1.5}, and~\thetag{1.6} we
identify the middle arrow, which is a pushforward in $K$-homology, as the
holomorphic induction
  $$ \aligned
      \ZZ[\alpha ,\alpha \inv ] &\longrightarrow \;\ZZ[\sigma _1]\;\times\;
     \ZZ[\sigma _2] \\
      \rho & \longmapsto \bigl(\ind(\rho ),\ind(\alpha ^m\rho ) \bigr)
     .\endaligned \tag{1.10} $$
The image of the generators of~$R(\TT)$ are expressed in terms of symmetric
products of the standard representation:
  $$ \aligned
      \alpha \inv &\longrightarrow \bigl(0,\Sym^{m-1}(\sigma _2) \bigr) \\ 
      1 &\longmapsto \bigl(1,\Sym^m(\sigma _2) \bigr).\endaligned 
     $$
Thus \thetag{1.10}~is injective and the cokernel is easy to compute, since
\thetag{1.10}~is a map of $R(SU_2)$-modules:
  $$ \aligned
      K^G_{1+\zeta }&=0, \\ 
      K^G_{0+\zeta }&\cong R(SU_2)/\langle \Sym^{m-1}(\sigma _2)
     \rangle.\endaligned \tag{1.11} $$
By Poincar\'e duality we obtain the twisted equivariant $K$-cohomology, which
is nonzero in odd degrees.  Since \thetag{1.9}~is a sequence of
$R(SU_2)$-module maps, \thetag{1.11}~is an isomorphism of rings: the twisted
equivariant $K$-theory is a quotient of the representation ring.  To make
contact with the theorem stated in the introduction, take $m=k+2$.  Then this
quotient of the representation ring is precisely the Verlinde algebra
of~$SU_2$ at level~$k$ (see~\cite{V}).

 \newpage
 \head
 \S{2} Chern-Simons Theory Revisited
 \endhead
 \comment
 lasteqno 2@  3
 \endcomment

A characteristic property of a topological quantum field theory in
$n$~dimensions is the {\it gluing law\/}.  The {\it partition function\/} is
a functorial assignment
  $$ X^n\longmapsto   Z(X)\quad \text{element of $\CC$} \tag{2.1} $$
of a complex number to a closed $n$-manifold~$X$.  Now if the manifold~$X$ is
cut by a closed hypersurface~$Y$ into a union of manifolds~$X_1,X_2$ with
boundary, then there are invariants~$Z(X_1),Z(X_2)$ so that $Z(X)=Z(X_1)\cdot
Z(X_2)$.  But $Z(X_i)$ is not a single complex number.  Rather, for each
``boundary condition''~$\alpha $ on~$Y$ we obtain complex
numbers~$Z(X_i)_\alpha $, and the invariant~$Z(X_i)$ is the vector of these
complex numbers.  (In {\it topological\/} quantum field theories there is
often a {\it finite\/} basis of~$\alpha $.)  The complex vector space of
these vectors is then functorially attached to the closed manifold~$Y$:
  $$ Y^{n-1}\longmapsto  E(Y)\quad \text{$\CC$-module}. \tag{2.2} $$
The standard story ends here.  But, as many people observed (see~\cite{F1},
for example) it is beneficial to go further and consider
hypersurfaces~$S\subset Y$ which express $Y$ as the union of submanifolds
$Y_1, Y_2$ with boundary.  Then the locality property leads us by analogy to
define a vector of complex vector spaces $E(Y_i) = \bigl(E(Y_i)\mstrut _\beta
\bigr)$ indexed by some boundary conditions~$\beta $ on~$S$.  It is natural
to view this as an element of a ``$K$-module'' functorially attached to~$S$,
which now has codimension~2 compared to the top dimension~$n$:
  $$ S^{n-2}\longmapsto  \sE(S)\quad \text{$K$-module}. \tag{2.3} $$
Here we view~$K$ heuristically as the category of finite dimensional
($\zt$-graded) complex vector spaces.
 
One can immediately see that a $K$-module is a category (whereas a
$\CC$-module is a set), but it has an additive structure as well.  The
additive structure on~$E(Y)$ in~\thetag{2.2} is crucial in quantum
mechanics---it encodes the superposition of states.  It stands to reason that
$\sE(S)$~should have an additive structure as well.  The discussion in this
section is completely heuristic, so we do not give a formal definition of
$K$-modules; one definition in the literature, where it is called a
``2-vector space,'' runs several pages.  Rather, we give an example.  Let
$F$~be a finite group.  Then the category~$\sE(F)$ of complex finite
dimensional representations of~$F$ is a $K$-module.  Addition is the direct
sum of representations, and scalar multiplication is tensor product with a
vector space on which $F$~acts trivially.  Since we can multiply
representations using the tensor product, $\sE(F)$~is a $K$-algebra.  Even
better, there is a sort of Hilbert space structure: the inner product of
representations~$E_1,E_2$ is the vector space of $F$-invariant maps $E_1\to
E_2$.
 
One can continue down this path, slicing and dicing manifolds of increasing
codimension.  This quickly leads to multicategories.\footnote{so perhaps also
to multiheadaches.  I define a mathematician's {\it category number\/} to be
the largest~$n$ such that he/she can think about $n$-categories without
getting a migraine.}  Such structures probably exist in any quantum field
theory, since locality is a characteristic property, but the infinite
dimensionality renders them of little use.
 
Now specialize to~$n=3$.  Since every closed 1-manifold is a finite union of
circles, there is only one interesting $K$-module~$\sE=\sE(S^1)$ in the
theory.  Elementary topology in one and two dimensions gives extra structure
to~$\sE$.  For example, in a theory of {\it oriented\/} manifolds reflection
in~$S^1$ gives an involution on~$\sE$.  The disk defines a distinguished
element of~$\sE$, and the pair of pants a multiplication on~$\sE$.  The
standard arguments which demonstrate that the $\CC$-module attached to~$S^1$
in a two-dimensional field theory is a {\it Frobenius algebra\/} over~$\CC$
lead us to the conclusion that in a three-dimensional theory $\sE$~is a
``Frobenius algebra over~$K$.''  We do not pretend to give a formal
definition, but note that in the literature the object attached to the circle
is usually called a ``modular tensor category,'' or some close variation.  In
fact, there is a theorem that the entire three-dimensional topological
quantum field theory can be reconstructed from this algebraic data
(see~\cite{T}).
 
In this context the {\it Verlinde algebra\/} is the Grothendieck group of the
$K$-algebra~$\sE$.  In particular, it is an algebra over the Grothendieck
group~$K^\bul(pt)$ of~$K$.
 
Three-dimensional Chern-Simons theory is the most well-studied
example.\footnote{There is a subtlety we ignore: In Chern-Simons theory the
category of oriented manifolds is centrally extended to oriented manifolds
with ``$p_1$-structure.''}  It is defined for any compact Lie group~$G$
starting with a ``level''~$k\in H^4(BG;\ZZ)$.  One gains insight from the
case when $G$~is finite, and this was much studied.  In particular, the
$K$-module~$\sE$ has a simple description in that case.  If $k=0$ then
$\sE=\Vect_G(G)$, the category of $G$-equivariant complex vector bundles
over~$G$.  The $K$-module structure uses pointwise direct sum and tensor
product; the algebra structure is induced from the multiplication $G\times
G\to G$ by pushforward.  The Grothendieck group is the equivariant $K$-theory
group~$K_G(G)$.  This is the simplest case of our main theorem, and it was
well-known ten years ago.  (The observation that equivariant $K$-theory
enters is credited to Lusztig.)  A similar picture holds if the level~$k$ is
nonzero.  Namely, there is a central extension of the category ``$G$~acting
on~$G$'' by a transgression of~$k$: to each pair~$(x,g)$ of elements of~$G$,
thought of as an arrow from~$x$ to~$gxg\inv $, is attached a complex line,
and these lines multiply suitably under composition of arrows.  The
$K$-module~$\sE$ consists of vector bundles over~$G$ suitably equivariant
under this central extension.  The precise statements appear in~\cite{F2},
but we failed to recognize the Grothendieck group as a twisted $K$-theory
group.  That belated realization led to the main theorem.
 
In many cases the partition function~\thetag{2.1} of a quantum field theory
is formally written as a functional integral in terms of classical fields and
a classical action.  One of the key ideas in~\cite{F2} is to extend the
notion of classical action in an $n$-dimensional field theory to manifolds of
dimension~$<n$, and to view the quantum Hilbert space~\thetag{2.2} and the
quantum $K$-module~\thetag{2.3} as defined by an extended notion of the
functional integral.  For gauge theories with finite structure group the
functional integral reduces to a finite sum, and so can be explicitly
computed.  This led us to the description of~$\sE$ above.  For any compact
group~$G$ we can best describe the extended {\it classical\/} action in terms
of {\it differential cohomology\/}.  This marriage of integral cohomology and
differential forms first appeared in differential geometry as {\it
Cheeger-Simons differential characters\/}, and also goes by the name of {\it
smooth Deligne cohomology\/}.  A treatment in the spirit needed here is given
in~\cite{HS}.  The level~$\lambda \in H^4(BG;\ZZ)$ lifts uniquely to a
differential cohomology class in~$\Hh^4(BG)$, and we fix a particular
cocycle~$\lh$ which represents it.  Then if $P\to M$ is any principal
$G$-bundle with connection~$A$, there is a degree~4 differential
cocycle~$\lh(A)$ on~$M$ which carries the information of a characteristic
class of~$P$ in~$H^4(M)$ and the Chern-Weil form of~$A$ in~$\Omega ^4(M)$.
If now $M\to T$~is an oriented fiber bundle of compact manifolds with a
$G$-connection~$A$ on~$M$, then we can integrate~$\lh(A)$ over the fibers
of~$M\to T$ to construct invariants.  For $X\to T$ with closed fibers of
dimension~3 that integral is a cocycle for~$\Hh^1(T)$, which is simply a map
$T\to \TT$, where $\TT$~is the circle group.  This is the exponentiated
Chern-Simons invariant.  It is the integrand of the functional integral which
defines the partition function~\thetag{2.1} of quantum Chern-Simons
theory~\cite{W}.  Note that in defining the functional integral (over the
universal family) we must extend the codomain of the Chern-Simons invariant
from~$\TT$ to~$\CC$.  Of course, that functional integral is a formal
expression; the only mathematically rigorous definition of these invariants
heretofore begin with the algebraic data in dimension one and build up from
there.
 
Next, consider a fiber bundle $Y\to T$ whose fibers are closed oriented
2-manifolds, and suppose we have a $G$-connection~$A$ on~$Y$.  Now the
integral of~$\lh(A)$ over the fibers is a cocycle for~$\Hh^2(T)$, which we
may view as a $\TT$-bundle with connection over~$T$.  It is the appropriate
``classical action'' in dimension~2 for this theory.  The associated
hermitian line bundle over the universal family of {\it flat\/} connections
is the ``prequantum line bundle'' in the theory of geometric quantization.
In that story one chooses a ``polarization'' of the symplectic manifold of
flat connections, then constructs a Hilbert space of compatible 1/2-forms
with coefficients in the prequantum line bundle.  The most natural
polarization in this case is complex, in which case we take the space of
holomorphic sections of the prequantum line bundle tensored with the square
root of an appropriate canonical line bundle.  It is conceptually useful to
regard the geometric quantization procedure as a functional integral of the
classical action in dimension~2. 
 
Finally, we pass to fiber bundles $S\to T$ whose fibers are oriented circles.
Then the appropriate classical action associated to a $G$-connection~$A$
on~$S$ is a cocycle for~$\Hh^3(T)$, which may be thought of as a
``$\TT$-gerbe with connection''.  To define the quantum
invariant~\thetag{2.3} we let $T$~be the universal family of connections
on~$S^1$ modulo gauge equivalence.  If we fix a basepoint on~$S^1$, then a
connection is determined up to isomorphism by its holonomy; changing the
basepoint conjugates the holonomy.  So the category of $G$-connections
on~$S^1$ is equivalent to the category which expresses the action of~$G$ on
itself by conjugation.  Now just as in the previous paragraph we passed from
a bundle of $\TT$-torsors to a bundle of $\CC$-modules, here we pass from a
bundle of $\TT$-gerbes to a bundle of $K$-modules.\footnote{Let~$\sT$ denote
the category of $\TT$-torsors.  A $\TT$-gerbe is a $\sT$-torsor, and now the
analogy of linearizations is better: in the first instance we replace the
group~$\TT$ by the ring~$\CC$; in the second we replace the group~$\sT$ by
the ring~$K$.}  Formally, then, we expect~$\sE$ to be a suitable space of
sections of this bundle and its Grothendieck group---the Verlinde
algebra---to be the space of homotopy classes of sections.  From the
exposition in~\S{1} we recognize it as a twisted equivariant $K$-theory
group.  We have almost arrived at the statement of the main theorem, but
there are two caveats.  First, the twisting class~$\zeta $ which appears is
deviates from the Chern-Simons action by a universal ``adjoint shift.''  That
shift is analogous to $1/2$-form twist in geometric quantization.  It is
intriguing to ponder if there is a story---even heuristic---about
polarizations in this context which leads to the adjoint shift.  Second, the
nonzero $K$-cohomology appears in degree~$\dim G$, whereas these ideas lead
one to suspect it would appear in degree~0.  We have yet to find an intuitive
explanation of this.
 
The mathematical proof of the main theorem~\cite{FHT} does not involve any of
these heuristics.

\enddocument